# Ensemble Differential Evolution with Simulation-Based Hybridization and Self-Adaptation for Inventory Management Under Uncertainty

Sarit Maitra[1][0000-1111-2222-3333], Vivek Mishra[2][1111-2222-3333-4444], Sukanya Kundu[3][2222-3333-4444-5555], Maitreyee Das[4][3333-4444-5555-6666]

[1,2,3,4] Alliance University, Bengaluru, India
sarit.maitra@gmail.com

**Abstract.** This study proposes an Ensemble Differential Evolution with Simulation-Based Hybridization and Self-Adaptation (EDESH-SA) approach for inventory management (IM) under uncertainty. In this study, DE with multiple runs is combined with a simulation-based hybridization method that includes a self-adaptive mechanism that dynamically alters mutation and crossover rates based on the success or failure of each iteration. Due to its adaptability, the algorithm is able to handle the complexity and uncertainty present in IM. Utilizing Monte Carlo Simulation (MCS), the continuous review (CR) inventory strategy is examined while accounting for stochasticity and various demand scenarios. This simulation-based approach enables a realistic assessment of the proposed algorithm's applicability in resolving the challenges faced by IM in practical settings. The empirical findings demonstrate the potential of the proposed method to improve the financial performance of IM and optimize large search spaces. The study makes use of performance testing with the Ackley function and Sensitivity Analysis with Perturbations to investigate how changes in variables affect the objective value. This analysis provides valuable insights into the behavior and robustness of the algorithm.

**Keywords:** Ackley Function, Differential Evolution, Evolutionary Algorithm, Ensemble Optimization, Inventory Management, Self-Adaptive;

## 1 Introduction

Inventory management (IM) encounters significant challenges in dealing with uncertainty and stochastic demand. Within the field of operations research, optimizing the IM process has emerged as a prominent area of study. Since the objective functions in IM are multidimensional, non-convex, non-differentiable, non-continuous, and have several local optima, they pose a tough challenge. Analytically solving these problems is not only challenging but, in many cases, impossible. Consequently, metaheuristic algorithms have been identified as advantageous over traditional algorithms in solving the optimal solutions to such problems (e.g., [15]; [17]; [3], etc.). However, they do not guarantee finding the globally optimal solution. This means that they may not always converge on the globally optimal solution. They rely on heuristics and search strategies that balance exploration and exploitation, aiming to improve the solution iteratively. Therefore, there is always a possibility of getting trapped in local optima, especially in complex and multimodal optimization problems like IM, where the suboptimal solutions are in the vicinity of the current search space.

We explore the scope of the existing DE approach to solving optimization problems. Through this work, we introduce EDESH-SA to overcome the trap of local minima and enhance the performance of the algorithm. DE has a competitive advantage over related EAs in terms of floating-point encoding [35]. According to a recent study, 158 out of 192 papers were published between 2016 and 2021, showing that academics have improved DE to increase its effectiveness and efficiency in handling a variety of optimization challenges [2]. This clearly indicates a consistent usage pattern of DE for optimization purposes.

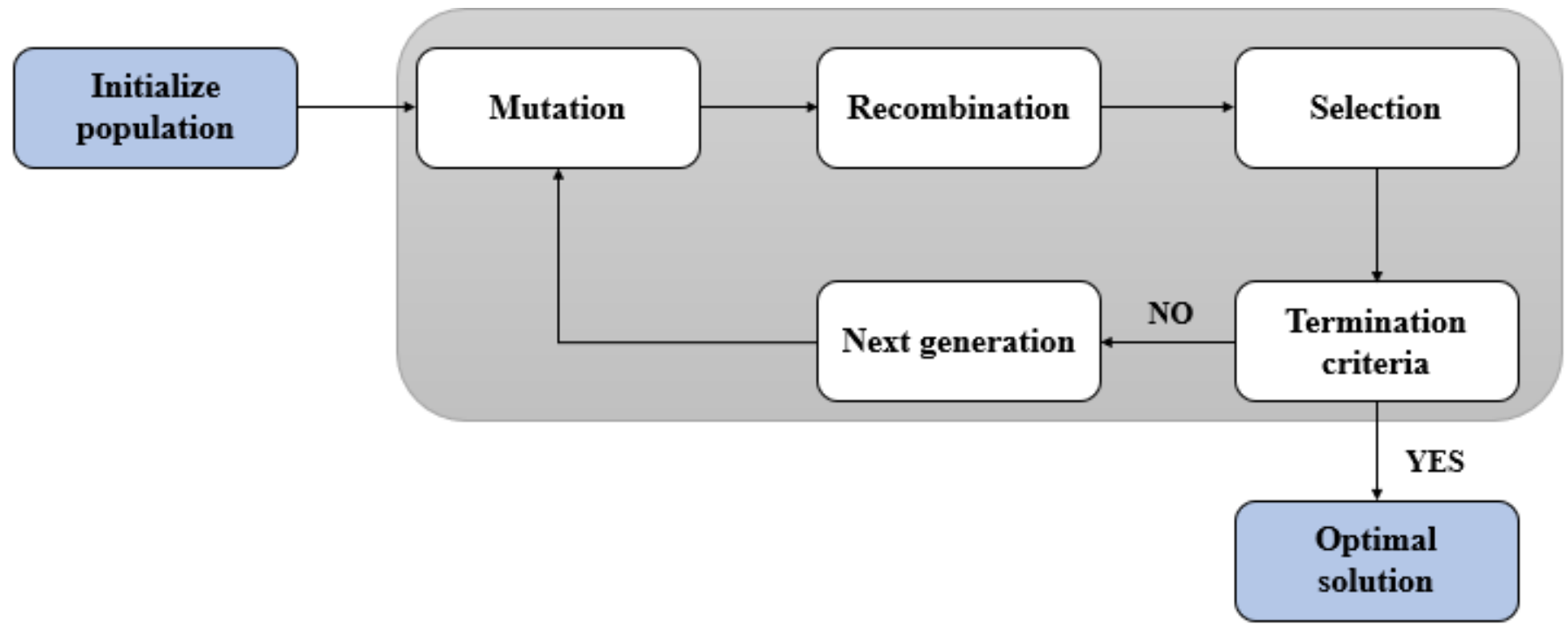

**Fig. 1.** Evolutionary algorithm procedure

At the center of the DE algorithm is the sequential optimization procedure as well as the interaction of numerous components, including selection, mutation, cross-over, and fitness evaluation (see Fig. 1). To optimize the decision variables (control factors), the proposed EDESH-SA integrates both control and noise elements. It strengthens the resilience and convergence of the optimization process while increasing the diversity of solutions. Ackley function and sensitivity analysis further establish the strength of the proposed solution.

The key contribution of this work lies in the development of EDESH-SA which combines statistical analysis, simulation, and optimization to overcome the shortcomings of the conventional DE algorithm. It provides an effective framework for decision-making, improves the performance of the DE algorithm, and offers valuable insights into the optimization process in the context of IM under uncertainty.

## 2 Previous Work

Optimization techniques provide a comprehensive framework for determining the most effective inventory policies, optimal inventory levels, and other key characteristics. The importance of optimization has been emphasized by several authors (e.g., [23]; [22]; [11]; [4]). Studies (e.g., [27]; [28]; [23]) and industry reports suggest that IM costs can range from ≥ 20–40% of the total SCM costs. Researchers have conducted a study on optimization under uncertainty and proposed a simulation-optimization approach [12]. Their findings presented an analytical report showing a 16% reduction in supply and managing costs by implementing the optimal policy. Probability distribution and MCS have been part of IM for a long time (e.g., [14]; [29]; etc.). Some of the recent work strengthens the argument by emphasizing the use of MCS in their work ([19] and [20]).

IM under uncertainty is challenging to solve due to the non-linearity of the model and several local optimum solutions ([9] and [16]). As a result, metaheuristic algorithms are frequently employed as a powerful solution for IM ([13]; [8]). They found that DE has received significant attention from researchers and practitioners over the last two decades [25]. Since the introduction of DE [31], it has undergone a plethora of performance enhancements to deal with complex optimization problems. However, there are still research gaps that need continued work on the advancement of DE.

It has been found that the ensemble-based technique helps DE choose its own parameter control, changing it to fit the most viable solution [32]. Ensemble-based DE could be considered a dynamic attempt to handle mutation and crossover techniques in a more constrained setting while maintaining a constant convergence rate. Though self-adaptive DE algorithms with multiple strategies were proposed [7], some authors emphasized the need for self-adaptive DE with minimum intervention [32]. It was noticed that, despite their infrequent explicit application, self-adaptation procedures are a novel adaptation tool with a high level of robustness [21].

While all the researchers acknowledged the challenges and opportunities in IM under uncertainty, the field is evolving, especially post pandemic. A growing body of research can be seen incorporating stochastic demand models, evaluating the impact of inventory policies, and exploring the capabilities of various optimization techniques. Metaheuristic optimization is a promising area of research (e.g., [1]; [10]; [33]; [8]; [29]; etc.).

IM is an evolving field where the goal is to enhance decision-making processes and address the challenges posed by uncertainties, such as those experienced during and post-pandemic. Though there is an existing body of knowledge on MCS, ensemble-based techniques, and self-adaptive DE algorithms separately, a combination of these approaches in one unified framework has not been explored. There is an opportunity to add value by proposing and investigating a novel approach that integrates all.

## 3 Methodology

Through the research gap, this methodology combines statistical analysis, MCS, ensemble, and self-adaptive DE. The data for demand information was collected for multiple products over the last 12 months. This was necessary to mathematically model the stochastic demand. We have calculated the best inventory levels that balance the costs of retaining stocks against the costs of stockouts by modelling numerous scenarios. These levels of inventory balance the probabilities of experiencing various levels of demand. It was assumed here that order sizes would have unknown values from a lognormal distribution, which is often the case in real-life scenarios. By running the simulations, a certain number of times (10,000 times in this case), the algorithm tries a variety of scenarios and fine-tunes the inventory policy to maximize profits. To model the inventory's performance over a year, we ran this simulation for a full 365 days.

The major sources of uncertainty are identified here:

- Unpredictable purchase behavior of the customer:

$\rho = (number\ of\ orders\ last\ year)/(number\ of\ working\ days)$, where $\rho$ is the probability of a customer placing an order on any given day. However, this model is not perfect, and there may be some level of error in the estimated probability of a customer placing an order.

- Variability in order size is modelled using a log-normal distribution and this is frequent practice in real-life scenarios. While past sales data is used to predict distribution parameters, actual order sizes may still differ greatly from these projections.

[26] conducted an extensive literature review of over seven decades and discovered that the CR is the most employed policy in stochastic inventory literature. Taking a clue from their work, we examined the inventory policy by inducing CR in our empirical analysis.

## 4 Data

The business case selected here examines the sale of four distinct products and considers adopting a suitable IM policy. The historical demand data is used to calculate the central tendency of the data, such as the mean demand during lead time, the standard deviation, and the minimum and maximum demand levels (see Table 1).

- PrA, PrB, PrC, and PrD are four distinct products.
- PC = purchase cost of one unit of the item from the supplier.
- LT = lead time to deliver the item after placing an order.
- Size = quantity of each item.
- SP = selling price of each item.
- SS = starting stock of each item in the inventory.
- $\mu$ = average demand for each item over a given period.
- σ = standard deviation of demand
- OC = order cost.
- HC = holding cost of one unit of inventory for a given period.
- $\rho$ = the probability of a stock-out event occurring, i.e., the probability of demand exceeding the available inventory level.

- $LT_{demand}$ = the demand during the lead time.

**Table 1.** Summary statistics

| | Pr A | Pr B | Pr C | Pr D |
|---|---|---|---|---|
| PC | € 12 | € 7 | € 6 | € 37 |
| LT | 9 | 6 | 15 | 12 |
| Size | 0.57 | 0.05 | 0.53 | 1.05 |
| SP | € 16.10 | € 8.60 | € 10.20 | € 68 |
| SS | 2750 | 22500 | 5200 | 1400 |
| σ | 37.32 | 26.45 | 31.08 | 3.21 |
| $\boldsymbol{\mu}$ | 103.50 | 648.55 | 201.68 | 150.06 |
| OC | € 1000 | € 1200 | € 1000 | € 1200 |
| HC | € 20 | € 20 | € 20 | € 20 |
| $\boldsymbol{\rho}$ | 0.76 | 1.00 | 0.70 | 0.23 |
| $LT_{demand}$ | 705 | 3891 | 2266 | 785 |

The shape of the KDE curves (see Fig. 2) provides insights into the underlying stochastic distribution of the data. The isolated peaks in the curves show potential outliers in the demand data.

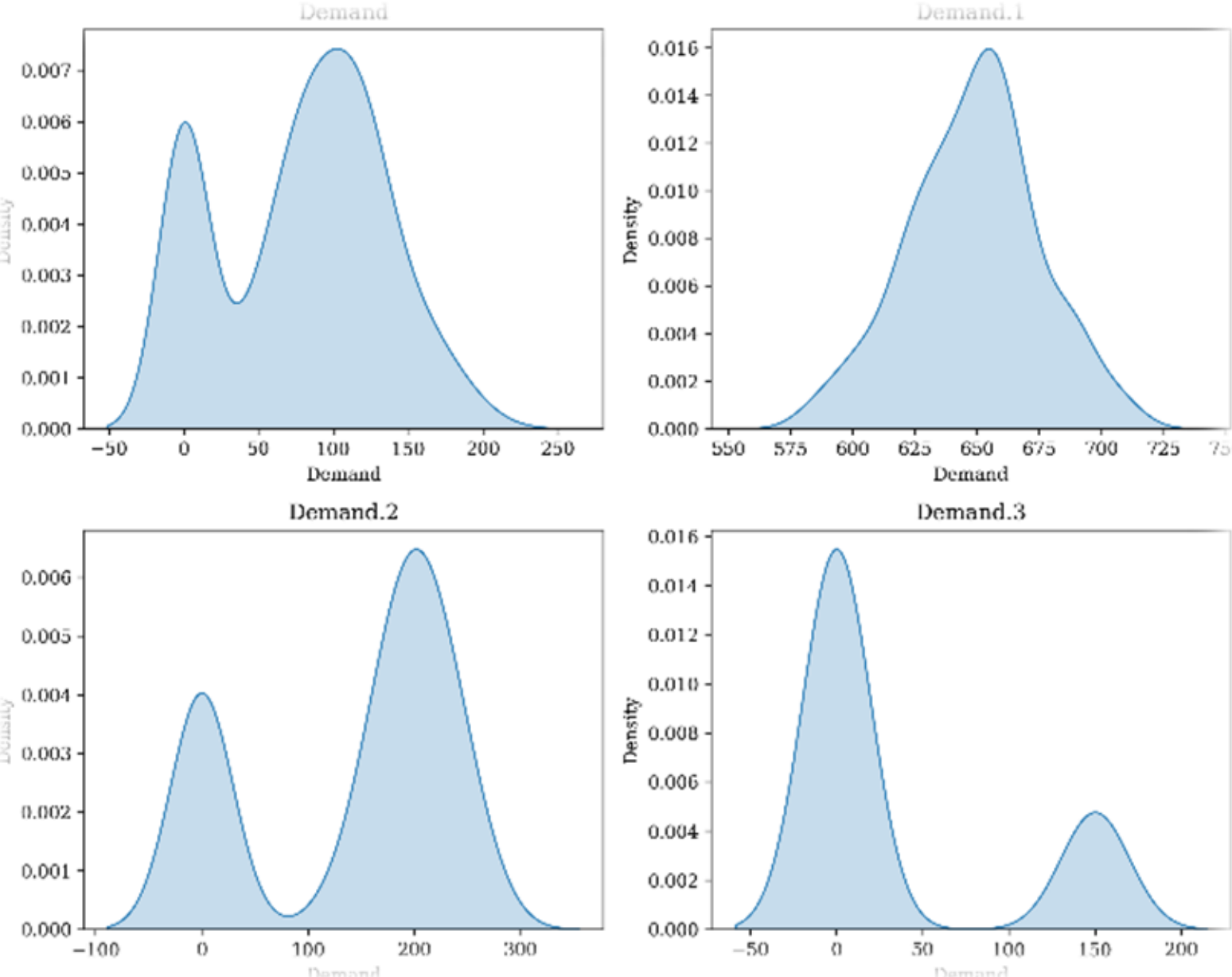

**Fig. 2.** KDE plots of demand distribution (365 days)

To demonstrate the methodology, we have used a contrived example that has been streamlined and simplified. We are able to concentrate on the essential components to make it simpler to comprehend and analyze the underlying concepts by removing superfluous complexity.

## 5 Ensemble Differential Evolution with Simulation-based Hybridization and Self-Adaptation

The proposed approach starts with finding the combination of control and noise factors the maximize the revenue and minimize the cost. The objective function in this approach calculates the profit as the weighted difference between the revenue and the cost. The DE algorithm is used to search for the optimal solution that achieves the highest profit. Fig. 3 displays the workflow of the proposed solution, which is an improvised version of DE.

The control factors in Fig. 3 include LT, Size, and SS. These are within the control of the decision-makers and can be adjusted to optimize the objective function. The noise factors introduce uncertainty into the optimization problem. These factors include PC and $LT_{demand}$, which are not directly under the control of the decision-makers and may vary, leading to variability in the optimization results. By incorporating both control and noise factors, the optimization problem accounts for both controllable and uncontrollable variables, enabling the determination of the optimal solution considering the uncertainties present in the system.

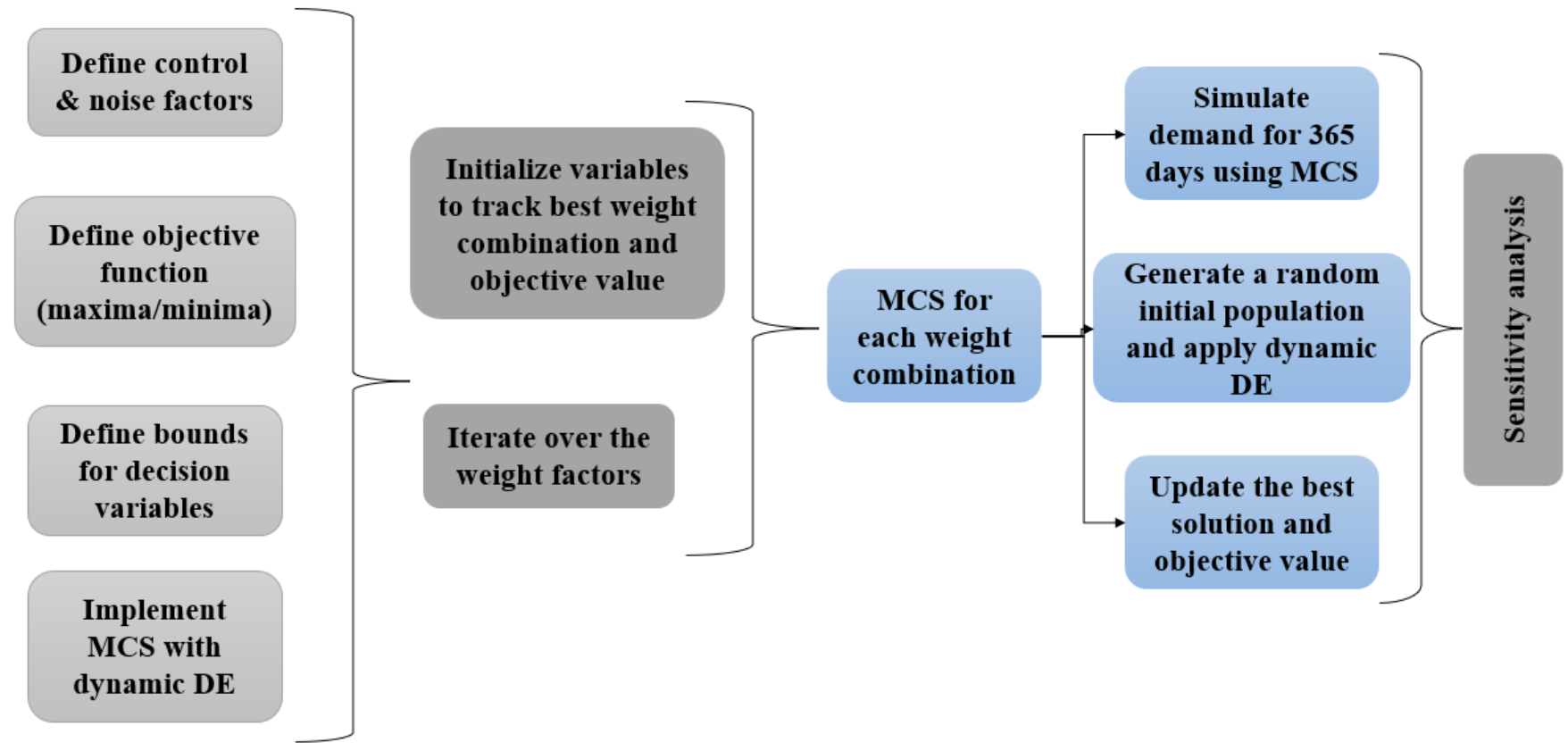

**Fig. 3.** Workflow of proposed solution

The reorder point, order quantity, and safety stock play crucial roles for CR. When the inventory level falls below the reorder point, instead of placing an order and waiting for delivery, the inventory is immediately replenished using the available starting stock and the lead time.

The reorder point (RP) is calculated as:

$$\text{RP} = LT_{demand} + \text{SSF} * \sigma \tag{1}$$

where, Safety stock factor (SSF) accounts for the variability in demand during the lead time and is estimated from the $\Sigma LT * \mu$

The order quantity (OQ):

$$\left(\frac{RP - SS}{OQ_{min}}\right) * OQ_{min}[\text{x }[0]] \tag{2}$$

where, x [0] represents the index of the control factor, $OQ_{min}$= Minimum order quantity corresponding to the control factor index x [0].

The revenue is calculated by multiplying the selling price of each product by the number of units sold for that product over the entire year. So, for each product, we add up the sales quantities for all 365 days and multiply them by the selling price. Annual profit is formulated as:

$$SP_i \sum_{t=1}^{365} S_{i,t} - \left\{\left(\frac{20HC_i}{365}\right) \sum_{t=1}^{365} I_{i,t} + N_i OC_{o,i} + \sum_{t=1}^{365} SC_i P_{i,t}\right\} \tag{3}$$

here,

- $SP_i$ = selling price of product i, $\sum_{t=1}^{365} S_{i,t}$ = sum of the sales quantities (S) for product i over the 365 days of the year. This component represents the total revenue generated by selling the product throughout the year.
- $HC_i$ = unit holding cost for product i, $\left(\frac{20HC_i}{365}\right) \sum_{t=1}^{365} I_{i,t}$ = holding cost for product i over the 365 days of the year; here it is 20 for all the products (Table 1). It considers the inventory levels (I) for each day.
- $N_i$ = number of orders placed for product i during the year.
- $OC_{(o,i)}$ = unit ordering cost for product i.
- $N_i OC_{o,i}$ calculates the ordering cost for product i based on the number of orders placed.
- $\sum_{t=1}^{365} SC_i P_{i,t}$ = stockout cost for product i over the 365 days of the year. It considers the inventory levels (P) and demand for each day.
- $SC_j$ = stockout cost for product i.

The net profit obtained from selling the product after accounting for various costs associated with inventory holding, ordering, and stockouts.

The total cost is computed based on the equation:

$$\text{cost} = \sum(\text{PC}_\text{i} * OQi) + \text{OC}_\text{i} * \text{I}(\text{OQ}_\text{i} > 0) + \sum\big(\text{HC}_\text{i} * \max(\text{x}_\text{i} - LT_{demand_\text{i}}, 0)\big) + \sum\big(\text{SC}_\text{i} * \max(LT_{demand_\text{i}} - \text{x}_\text{i}, 0)\big) \tag{4}$$

$SC = 20\% * \mu = [20.70, 129.71, 40.34, 30.01]$, assuming 20% is the stock out cost.

By performing MCS multiple times and averaging the costs and revenues, the objective function provides an estimation of the average cost and average revenue.

To find the inventory levels for each product that maximize the profit, we create an objective function that combines this revenue (Eq. 3) and cost components (Eq. 4). To balance maximizing revenue and minimizing costs, we introduce two weights (ω1 and ω2) that represent their relative importance.

$$f(x) \;=\; \omega_1 \;*\left[SP_i \sum_{t=1}^{365} S_{i,t} - \left\{\left(\frac{20V_i}{365}\right)\sum_{t=1}^{365} I_{i,t} +\; N_i C_{o,i} +\; \sum_{t=1}^{365} c_i P_{i,t}\right\}\right] - (1-\omega_1) * Cost_{optimal} \qquad (5)$$

Depending on the business priorities, the weight can be adjusted to prioritize revenue maximization or cost minimization. cost optimization or revenue maximization,

In this case, by incorporating both the sales revenue and average cost terms in the objective function ($f(x)$), the optimization algorithm aims to find the inventory levels x that not only minimize the average cost but also maximize the sales revenue. This provides a comprehensive approach to finding a better solution that considers both financial aspects (costs) and business performance (revenue).

### 5.1 Improvised Differential evolution

Centeno-Telleria et al.'s (2021) findings on the flexibility of DE, the benefits of adaptive techniques, and the emphasis on exploring the solution space are consistent with the characteristics of our work. The population matrix with the parameter vectors has the form:

$$x_{n,i}^{g} = \left\{x_{n,1}^{g}, x_{n,2}^{g}, x_{n,3}^{g}, \dots, x_{n,d}^{g}\right\} \qquad (6)$$

Where n = population size, g = generation and n = 1, 2, 3, …, N

$$x_{n,i} \;=\; x_{n,i}^{L} \;+\; \text{rand}\,() * \left(x_{n,i}^{u} - x_{n,i}^{l}\right),\;\; \text{i} \;=\; 1, 2, 3, \dots \text{ D and n} \;=\; 1, 2, 3, \dots \text{ N} \qquad (7)$$

- Three other vectors can be selected randomly from each parameter vector during the mutation phase such as $x_{r1n}^{g}$, $x_{r2n}^{g}$, and $x_{r3n}^{g}$. If we add the weighted difference of two of the vectors to the third, the equation becomes:

$$v_n^{g+1} \;=\; x_{r1n}^{g} \;+ \omega\left(x_{r2n}^{g} - x_{r3n}^{g}\right),\;\; \text{n} \;=\; 1, 2, 3, \dots \text{N} \qquad (8)$$

Where, $v_n^{g+1}$ = donor vector, $\omega$ is between 0 and 2.

*Instead, random selection, we employed MCS to choose the vectors based on their fitness values. The selection process is modified by selecting three vectors using MCS based on their fitness values: $x_{r1n}^{g}, x_{r2n}^{g}$, and $x_{r3n}^{g}$ and Eq (8) is modified as below:*

$$v_n^{g+1} \;=\; x_{r1n}^{g} \;+ \omega_n(g)\left(x_{r2n}^{g} - x_{r3n}^{g}\right), n \;=\; 1, 2, 3, \dots N \qquad (9)$$

*where $\omega_n(g)$ is the mutation factor for individual n at generation g. Instead of using a fixed ω value, we used a strategy to dynamically adapt ω based on the performance of the algorithm.*

- During recombination phase, a trial vector $u_{n,i}^{g+1}$ is developed from the target vector $\left(x_{n,i}^{g}\right)$ and the donor vector $\left(v_{n,i}^{g+1}\right)$

$$u_{n,i}^{g+1} = \begin{cases} v_{n,i}^{g+1} & if\ rand() \leq C_p\ or\ i = I_{rand} \quad i = 1,2,3,\dots D\ and \\ x_{n,i}^{g} & if\ rand() > C_p\ and\ i \neq I_{rand} \quad n = 1,2,3,\dots N \end{cases} \tag{10}$$

Where, $I_{rand}$ is an integer random number between [1, D] and $C_p$ is the recombination probability.

*We incorporate self-adaptive mechanisms, which allow the algorithm to dynamically adjust its parameters during the optimization process. Our proposed approach uses Adaptive Mutation Control, where we assign different mutation factors to different individuals based on their performance. This is done by associating a scaling factor with each vector and updating it based on their fitness improvement. We also employed Adaptive Recombination Control, where the recombination probability ($C_p$) is adaptively adjusted based on the current population's behavior. High recombination probabilities encourage exploration, while low values promote exploitation.*

- During the selection phase, the target vector $x_{n,i}^{g}$ is compared with the trial vector $u_{n,i}^{g+1}$ and the one with the lowest function value is chosen for the following generation.

$$x_{n}^{g+1} = \begin{cases} u_{n,i}^{g+1} & if\ f\left(u_{n}^{g+1}\right) < f\left(x_{n}^{g}\right) \\ x_{n}^{g} & otherwise \end{cases} \quad \mathrm{n} = 1,2,\dots \mathrm{N} \tag{11}$$

All the evolution phases are repeated until a predefined termination criterion is satisfied.

*Ensemble methods involve combining multiple instances of the optimization algorithm to improve overall performance. Each instance operates independently but shares information periodically to guide the evolution. By combining results from multiple (5) optimizers, the algorithm can avoid getting stuck in local minima and provide a better global optimal solution. The ensemble methods are integrated into DE:*

The improved version outlined above addresses several shortcomings of conventional DE algorithms, including convergence to local optima, lack of adaptability, limited exploration and exploitation capabilities, sensitivity to control parameters, and lack of robustness.

Table 2 displays a comparative analysis with conventional DE, wherein we have removed the MCS optimization loop and the adaptive mutation/crossover rate adjustments.

**Table 2.** Comparison of Optimization Results

| Description | Conventional DE | Proposed DE |
|---|---|---|
| Best weight combination | [0.7, 0.3] | [0.7, 0.3] |
| Best objective value | 474,230 | 58,760,547 |
| Mutation & cross-over | - | [0.5, 0.9] |
| Best solution | [0.93, 1.56, 0.98, 2.73, 2.87] | [3.0, 0.59, 0.42, 2.29, 1.89] |
| Annual profit | 30,335,547.85 | 67,102,330.00 |

The conventional DE achieved a best objective value of 474,230, while the proposed DE algorithm achieved a significantly higher best objective value of 58,760,547. The objective value represents the optimized balance between maximizing revenue and minimizing costs. Moreover, the conventional DE achieved an annual profit of 30,335,547.85, while the proposed DE algorithm achieved a higher annual profit of 67,102,330.00.

## 5.2 Algorithmic performance

Ackley function tests and sensitivity analyses were performed to gain insights into the performance characteristics of the differential evolution algorithm. These analyses provide valuable information about its efficiency, convergence behavior, and solution quality.

**Ackley Function**

. The Ackley function performance test shows the optimization algorithm's capacity to locate the best solution in a challenging and multidimensional search space. The contour plot in Fig. 4 visualizes the contours and valleys of the function, providing insights into the behavior of the objective function. The decision variables [2.63, 2.63, 1.00, 1.78, 1.33] correspond to the choice variable value that the DE algorithm finds to be the best. These values represent the optimized combination of control and noise factors that maximize the objective function. The objective value achieved is approximately -10.39, indicating the maximum value of the function obtained after optimizing the decision variables.

.

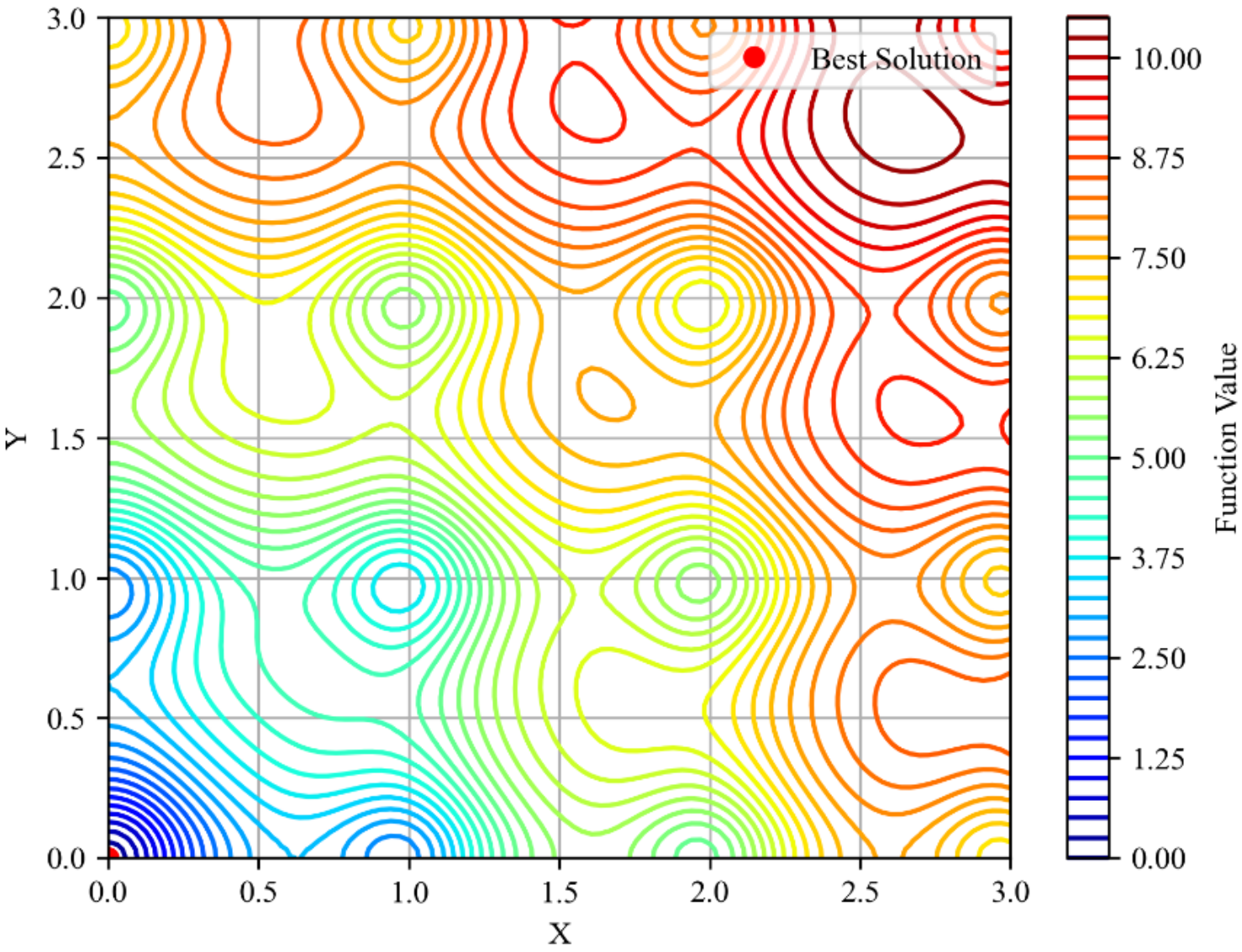

**Fig. 4.** Ackley Function Contour Plot

The 3D mesh grid plot (see Fig. 5) further illustrates the optimized combination of factors and the location of the best solution in the solution space. The red dot on the plot represents the optimal point where the function reaches its maximum value. The negative value of the objective function indicates the algorithm's success in finding a point of high fitness.

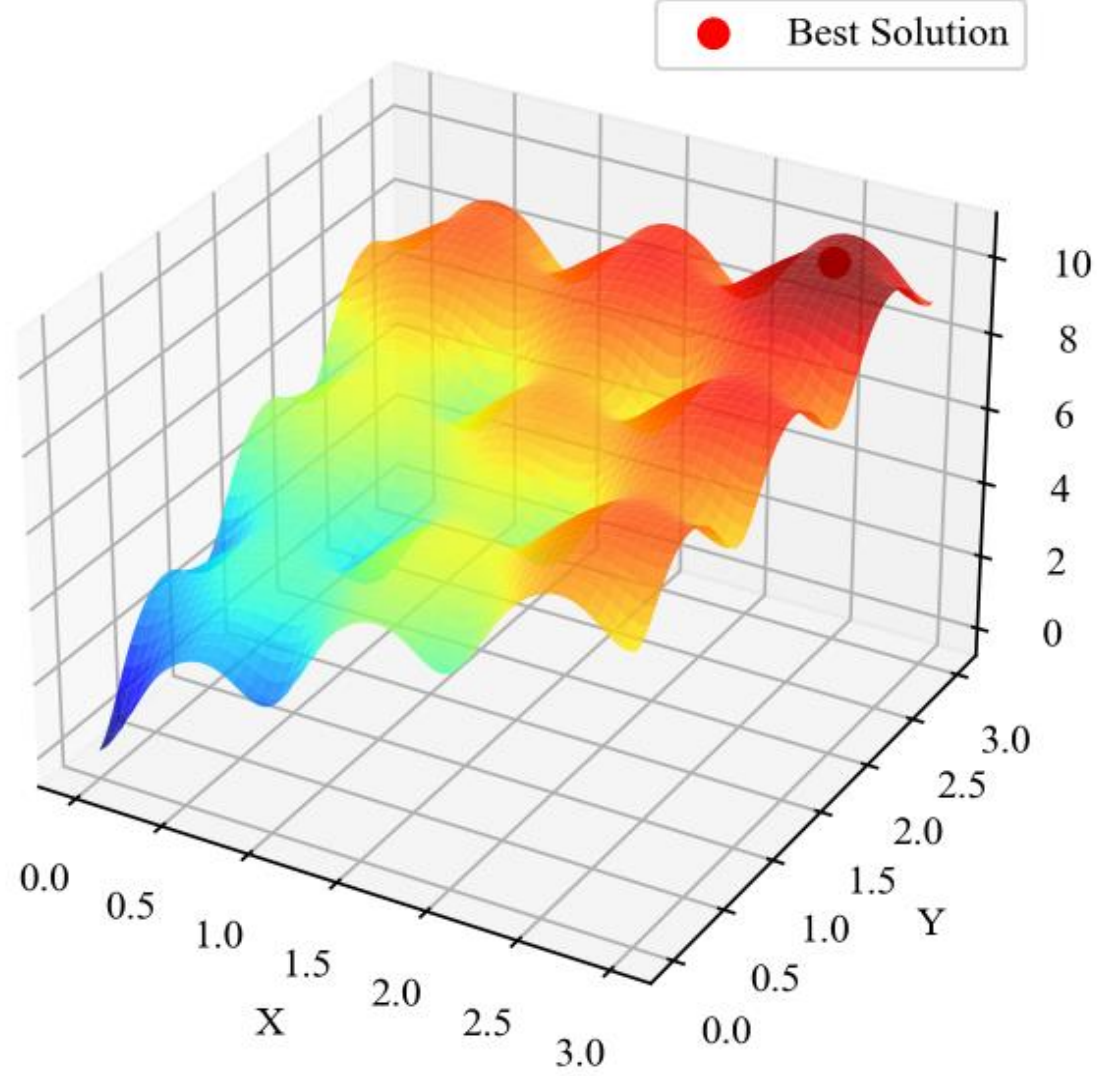

**Fig. 5.** Ackley Function 3D Mesh grid Plot

The Ackley function is known for its complex landscape with multiple local minima. Despite this challenge, the differential evolution algorithm effectively navigated the search space and found a point that maximizes the objective function. This demonstrates the algorithm's ability to handle complex and non-linear optimization problems.

**Sensitivity analysis**

. The use of perturbation analysis to calculate parameter sensitivities is an efficient method. Some researchers have made significant contributions to the field of sensitivity analysis and perturbation approaches in mathematical programming and optimization problems ([5], [34], and [18]). Their contributions provide a larger background and framework for the methodology used in our study. We compared the objective values before and after each variable perturbation to see how it affected the objective value. Here, the objective value worsens after the multiple perturbations, showing that changing the variable has a negative impact on the objective.

The x-axis (see Fig. 6) represents the variables being perturbed, and the y-axis represents the change in the objective value due to the perturbations. The blue line represents the change in the objective value when each variable is perturbed in the positive direction (increased). The orange line represents the change in the objective value when each variable is perturbed in the negative direction (decreased). If the blue line is below zero ($y = 0$) for a particular variable, it indicates that increasing the value of that variable leads to a worse objective value. Here, the orange line is below zero ($y = 0$) for a particular variable, which indicates that decreasing the value of that variable leads to a worse objective value.

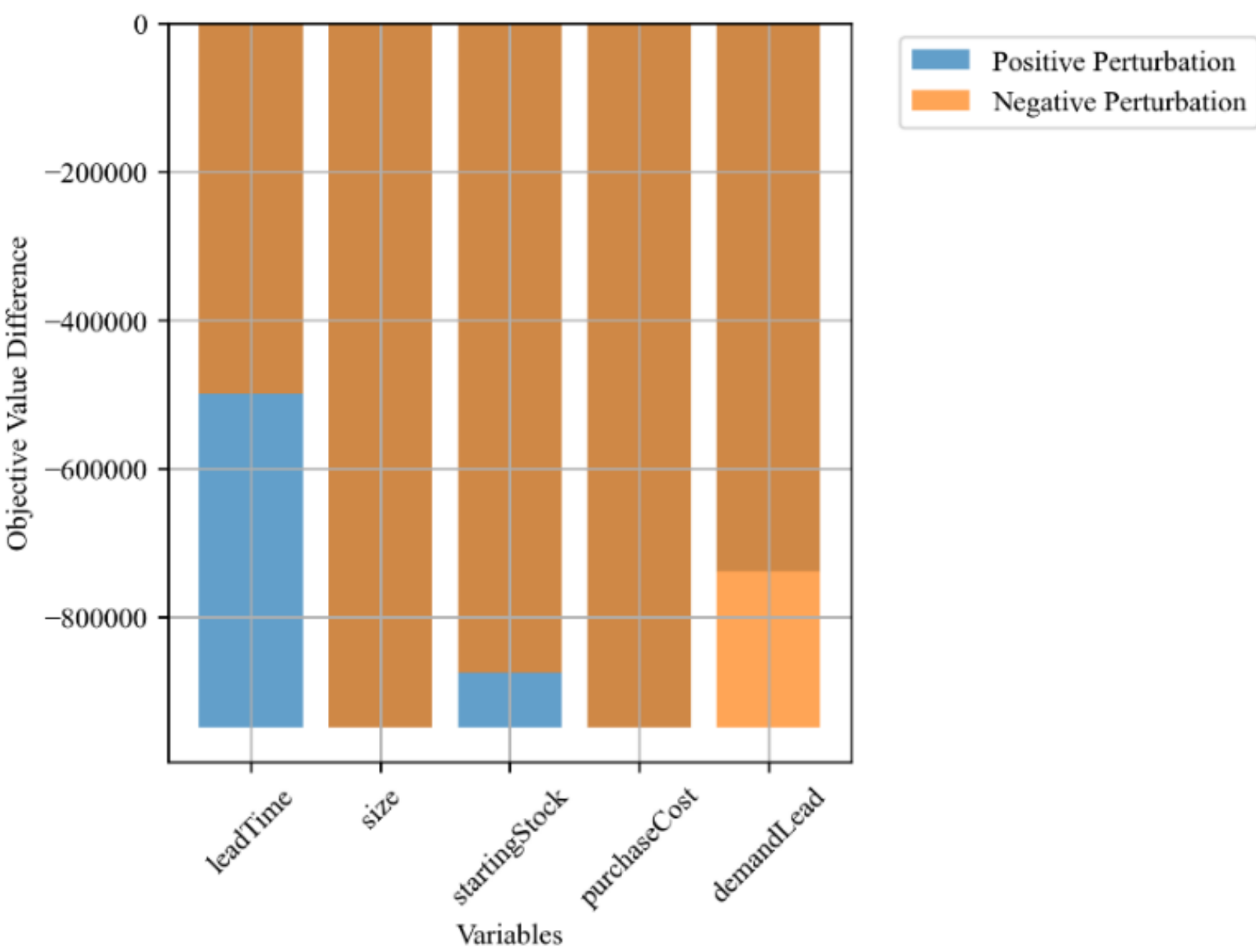

**Fig. 6.** Impact of Perturbations on Objective Value

Through this graphical representation, the sensitivity and impact of variable perturbations on the objective value are effectively demonstrated.

### 5.3 Limitations

The proposed approach primarily focuses on optimizing internal IM decisions, such as reorder points, order quantities, and safety stock levels. However, there are external factors that impact IM, such as supply chain disruptions, supplier reliability, transportation constraints, or changing market conditions. Incorporating these external factors into the optimization model could provide a more realistic and robust solution. Moreover, the provided solution does not address scalability, which can be achieved by adjusting parameters, employing parallelization techniques, optimizing algorithmic efficiency, and considering problem decomposition strategies. These factors should be considered when applying the solution to larger and more complex IM systems.

## 6 Conclusion

To improve the performance of the optimization technique for IM, this work provided a simulation-based hybridization strategy mixed with ensemble DE. By incorporating a self-adaptive approach, the algorithm can dynamically adjust mutation and crossover rates based on the success or failure of each iteration, further improving its effectiveness. The proposed approach considers both control and noise factors, optimizing decision variables while accounting for system stochasticity. Through MCS, the algorithm evaluates the performance of implemented inventory policies under different demand scenarios, providing a more comprehensive optimization process. The effectiveness of the proposed solution is demonstrated using an empirical analysis of the optimal combination of control and noise factors that maximize profit. The algorithm successfully navigates the complex search space of the Ackley function, finding the best solution despite the presence of multiple local minima. Furthermore, sensitivity analysis highlights the impact of perturbations on the objective value, providing insights into the sensitivity of the solution to changes in variables. This research advances optimization methods for IM with stochastic demands. By combining metaheuristic optimization, simulation, and sensitivity analysis, the proposed approach offers a practical and robust solution for real-world IM problems. Future research can further explore the application of these techniques in different business environments and extend the analysis to include additional factors and constraints.

## References

1. Abdi, A., Abdi, A., Fathollahi-Fard, A. M., & Hajiaghaei-Keshteli, M.: A set of calibrated metaheuristics to address a closed-loop supply chain network design problem under uncertainty. International Journal of Systems Science: Operations & Logistics, 8(1), 23-40 (2021).

2. Ahmad, M. F., Isa, N. A. M., Lim, W. H., & Ang, K. M.: Differential evolution: A recent review based on state-of-the-art works. Alexandria Engineering Journal, 61(5), 3831-3872 (2022).
3. Alejo-Reyes, A., Mendoza, A., & Olivares-Benitez, E.: Inventory replenishment decision model for the supplier selection problem using metaheuristic algorithms. 15, 1509-1535 (2021).
4. Bag, S., Wood, L. C., Mangla, S. K., & Luthra, S.: Procurement 4.0 and its implications on business process performance in a circular economy. Resources, conservation and recycling, 152, 104502 (2020).
5. Castillo, E., Conejo, A. J., Castillo, C., Mínguez, R., & Ortigosa, D.: Perturbation approach to sensitivity analysis in mathematical programming. Journal of Optimization Theory and Applications, 128, 49-74 (2006).
6. Centeno-Telleria, M., Zulueta, E., Fernandez-Gamiz, U., Teso-Fz-Betoño, D., & Teso-Fz-Betoño, A.: Differential evolution optimal parameters tuning with artificial neural network. Mathematics, 9(4), 427 (2021).
7. Deng, G., & Gu, X.: A hybrid discrete differential evolution algorithm for the no-idle permutation flow shop scheduling problem with makespan criterion. Computers & Operations Research, 39(9), 2152-2160 (2012).
8. Fahimnia, B., Davarzani, H., & Eshragh, A.: Planning of complex supply chains: A performance comparison of three meta-heuristic algorithms. Computers & Operations Research, 89, 241-252 (2018).
9. Fallahi, A., Bani, E. A., & Niaki, S. T. A.: A constrained multi-item EOQ inventory model for reusable items: Reinforcement learning-based differential evolution and particle swarm optimization. Expert Systems with Applications, 207, 118018 (2022).
10. Faramarzi-Oghani, S., Dolati Neghabadi, P., Talbi, E. G., & Tavakkoli-Moghaddam, R.: Meta-heuristics for sustainable supply chain management: A review. International Journal of Production Research, 1-31 (2022).
11. Fonseca, L. M., & Azevedo, A. L.: COVID-19: outcomes for global supply chains. Management & Marketing. Challenges for the Knowledge Society, 15(s1), 424-438 (2020).
12. Franco, C., & Alfonso-Lizarazo, E.: Optimization under uncertainty of the pharmaceutical supply chain in hospitals. Computers & Chemical Engineering, 135, 106689 (2020).
13. Goodarzian, F., Wamba, S. F., Mathiyazhagan, K., & Taghipour, A.: A new bi-objective green medicine supply chain network design under fuzzy environment: Hybrid metaheuristic algorithms. Computers & Industrial Engineering, 160, 107535 (2021)
14. Huijbregts, M. A.: Application of uncertainty and variability in LCA. The International Journal of Life Cycle Assessment, 3, 273-280 (1998).
15. Kang, H., Liu, R., Yao, Y., & Yu, F.: Improved Harris hawks optimization for non-convex function optimization and design optimization problems. Mathematics and Computers in Simulation, 204, 619-639 (2023).
16. Khalilpourazari, S., Pasandideh, S. H. R., & Niaki, S. T. A.: Optimization of multi-product economic production quantity model with partial backordering and physical constraints: SQP, SFS, SA, and WCA. Applied Soft Computing, 49, 770-791 (2016).
17. Khan, A. T., Cao, X., Brajevic, I., Stanimirovic, P. S., Katsikis, V. N., & Li, S.: Non-linear activated beetle antennae search: A novel technique for non-convex tax-aware portfolio optimization problem. Expert Systems with Applications, 197, 116631 (2022).
18. Kiran, R., Li, L., & Khandelwal, K.: Complex perturbation method for sensitivity analysis of nonlinear trusses. Journal of Structural Engineering, 143(1), 04016154 (2017).

19. Krumscheid, S., Nobile, F., & Pisaroni, M.: Quantifying uncertain system outputs via the multilevel Monte Carlo method—Part I: Central moment estimation. Journal of Computational Physics, 414, 109466 (2020).
20. Luengo, D., Martino, L., Bugallo, M., Elvira, V., & Särkkä, S.: A survey of Monte Carlo methods for parameter estimation. EURASIP Journal on Advances in Signal Processing, 2020(1), 1-62 (2020).
21. Meyer-Nieberg, S., & Beyer, H. G.: Self-adaptation in evolutionary algorithms. In Parameter setting in evolutionary algorithms. Berlin, Heidelberg: Springer Berlin Heidelberg, 47-75 (2007)
22. Modgil, S., Singh, R. K., & Hannibal, C.: Artificial intelligence for supply chain resilience: learning from Covid-19. The International Journal of Logistics Management, 33(4), 1246-1268 (2022).
23. Moons, K., Waeyenbergh, G., & Pintelon, L.: Measuring the logistics performance of internal hospital supply chains–a literature study. Omega, 82, 205-217 (2019).
24. Muller, M.: Essentials of inventory management. HarperCollins Leadership (2019).
25. Pant, M., Zaheer, H., Garcia-Hernandez, L., & Abraham, A.: Differential Evolution: A review of more than two decades of research. Engineering Applications of Artificial Intelligence, 90, 103479 (2020).
26. Perera, S. C., & Sethi, S. P.: A survey of stochastic inventory models with fixed costs: Optimality of (s, S) and (s, S)-type policies—Continuous-time case. Production and Operations Management, 32(1), 154-169 (2023).
27. Simchi-Levi, D., Kaminsky P., Simchi-Levi E.: Designing and managing the supply chain: concepts, strategies, and case studies, McGraw Hill Professional (2033).
28. Singh, D., & Verma, A.: Inventory management in supply chain. Materials Today: Proceedings, 5(2), 3867-3872 (2018).
29. Soleimani, H., & Kannan, G.: A hybrid particle swarm optimization and genetic algorithm for closed-loop supply chain network design in large-scale networks. Applied Mathematical Modelling, 39(14), 3990-4012 (2015).
30. Sonnemann, G. W., Schuhmacher, M., & Castells, F.: Uncertainty assessment by a Monte Carlo simulation in a life cycle inventory of electricity produced by a waste incinerator. Journal of cleaner production, 11(3), 279-292 (2003).
31. Storn, R., & Price, K. (1997). Differential evolution-a simple and efficient heuristic for global optimization over continuous spaces. Journal of global optimization, 11(4), 341.
32. Wang, S. L., Ng, T. F., & Morsidi, F.: Self-adaptive ensemble based differential evolution. International Journal of Machine Learning and Computing, 8(3), 286-293 (2018).
33. Wang, S., Wang, L., & Pi, Y.: A hybrid differential evolution algorithm for a stochastic location-inventory-delivery problem with joint replenishment. Data Science and Management, 5(3), 124-136 (2022).
34. Xu, G., & Burer, S.: Robust sensitivity analysis of the optimal value of linear programming. Optimization Methods and Software, 32(6), 1187-1205 (2017).
35. Zamuda, A., & Brest, J.: Self-adaptive control parameters׳ randomization frequency and propagations in differential evolution. Swarm and Evolutionary Computation, 25, 72-99 (2015).